\documentclass[journal]{IEEEtran}

\usepackage{bm}
\usepackage{wrapfig}
\usepackage{float}
\usepackage{amsmath}
\usepackage{amssymb}
\usepackage{multirow}
\usepackage{amsfonts}
\usepackage{mathrsfs}
\usepackage{graphicx}
\usepackage{epstopdf}
\usepackage{subcaption}
\usepackage{adjustbox}
\usepackage{textcomp}
\usepackage{autobreak}
\usepackage{tabularx}
\usepackage{comment}
\usepackage{caption}
\usepackage{lettrine}
\usepackage{cite}
\usepackage{setspace}
\usepackage{amsmath}
\usepackage{titlesec}
\usepackage{lipsum}
\usepackage{glossaries}
\usepackage[utf8]{inputenc}
\usepackage{tikz}
\usepackage{import}
\usepackage{dirtytalk}

\usepackage{hyperref}

\usepackage[ruled, lined, linesnumbered, commentsnumbered, longend]{algorithm2e}
\RestyleAlgo{ruled}
\usetikzlibrary{shapes.geometric, arrows}
\usetikzlibrary{fit}
\tikzstyle{arrow} = [thick,->,>=stealth]

\usepackage[utf8]{inputenc}
\usepackage{tikz}
\usepackage{import}
\usetikzlibrary{shapes.geometric, arrows}
\tikzstyle{startstop} = [rectangle, rounded corners, minimum width=3cm, minimum height=1cm,text centered, draw=black]
\tikzstyle{process} = [rectangle, minimum width=3cm, minimum height=1cm, text centered, text width=3cm, draw=black]
\tikzstyle{prow} = [rectangle, minimum width=5cm, minimum height=1cm, text centered, text width=5cm, draw=black]
\tikzstyle{prow2} = [rectangle, minimum width=8cm, minimum height=1cm, text centered, text width=8cm, draw=black]
\tikzstyle{decision} = [diamond, minimum width=3cm, minimum height=1cm, text centered, draw=black]
\tikzstyle{arrow} = [thick,->,>=stealth]

\titlespacing*{\section}
{0pt}{3.5ex plus 1ex minus .2ex}{2.0ex plus .2ex}

\ifCLASSINFOpdf
\else
\fi
\begin{document}

\title{Tight Data-Driven Linear Relaxations for Constraint Screening in Robust Unit Commitment}

\author{Mohamed~Awadalla and Fran\c{c}ois~Bouffard~\IEEEmembership{Senior Member,~IEEE}%
\thanks{This work was supported in part by IVADO, Montreal, QC, the Natural Sciences and Engineering Research Council of Canada, Ottawa, ON, and the Trottier Institute for Sustainability in Engineering and Design, Montreal, QC.}
\thanks{M. Awadalla and F. Bouffard are with the Department of Electrical and Computer Engineering, McGill University, Montreal, QC H3A~0E9, Canada and with the Groupe d'\'{e}tudes et de recherche en analyse des d\'{e}cisions (GERAD), Montreal, QC  H3T~1J4, Canada (emails: mohamed.awadalla@mail.mcgill.ca; francois.bouffard@mcgill.ca).}
}

\maketitle
\begin{abstract}
The daily operation of real-world power systems and their underlying markets relies on the timely solution of the unit commitment problem. However, given its computational complexity, several optimization-based methods have been proposed to lighten its problem formulation by removing redundant line flow constraints. These approaches often ignore the spatial couplings of renewable generation and demand, which have an inherent impact of market outcomes. Moreover, the elimination procedures primarily focus on the feasible region and exclude how the problem's objective function plays a role here. To address these pitfalls, we move to rule out \textit{redundant} and \textit{inactive} constraints over a tight linear programming relaxation of the original unit commitment feasibility region by adding valid inequality constraints. We extend the optimization-based approach called \textbf{\textit{umbrella constraint discovery}} through the enforcement of a consistency logic on the set of constraints by adding the proposed inequality constraints to the formulation. Hence, we reduce the conservativeness of the screening approach using the available historical data and thus lead to a tighter unit commitment formulation. Numerical tests are performed on standard IEEE test networks to substantiate the effectiveness of the proposed approach.


\end{abstract}

%
\IEEEpeerreviewmaketitle

\section*{Nomenclature}

The main symbols used in the paper are listed here. Other symbols will be defined as required.

\subsection{Sets and Indices}
\begin{IEEEdescription}[\IEEEusemathlabelsep\IEEEsetlabelwidth{$1,2,3,4$}]
\item[$\mathcal{N}$]{Set of buses, indexed by $n$ and of size $N$.}
\item[$\mathcal{L}$]{Set of transmission lines, indexed by $l$ and of size $L$.}
\item[$\mathcal{M}$]{Set of generating units, indexed by $m$ and of size $M$.}
\item[$\mathcal{M}_{n}$]{Set of generating units connected to bus $n$.}
\item[$\mathcal{T}$]{ Set of time periods, indexed by $t$ and of size $T$.}
\end{IEEEdescription}

\subsection{Variables}
\begin{IEEEdescription}[\IEEEusemathlabelsep\IEEEsetlabelwidth{$1,2,3,4$}]
\item[${g}_{m}$]{Power output dispatch of generating unit $m$.}
\item[${q}_{n}$]{Net injected power at node $n$.}
\item[${u}_{m}$]{Commitment of generating unit $m$.}
\end{IEEEdescription}

\subsection{Parameters}
\begin{IEEEdescription}[\IEEEusemathlabelsep\IEEEsetlabelwidth{$1,2,3,4$}]
\item[$c_{m}$]{Incremental production cost of generating unit $m$.}
\item[$d_{n}$]{Residual demand at node $n$.}
\item[$f_l^{\max}$]{Maximum flow capacity of transmission line $l$.} 
\item[$g_m^{\max}$]{Maximum power limit of generator $m$.}
\item[$g_m^{\min}$]{Minimum power limit of generator $m$.}
\item[$h_{l n}$]{Power transfer distribution factor (PTDF) of line $l$ for power injections at bus $n$.}
\end{IEEEdescription}


\section{Introduction}
\IEEEPARstart{D}{espite} the liberalization of the electricity sector, unit commitment (UC) is still a fundamental optimization tool for all major system operators in the United States and Canada for resolving market clearing underlying the daily operation planning of their respective power systems \cite{Sun2017, Zheng2015}. The solution to the UC problem determines the most economical operating schedule, given by the on/off commitment status and production levels of the generating units based on submitted offers to generate and bids to consume. The goal of the UC problem is to minimize system operational cost, while satisfying generation, network, and other security constraints. \par

Mathematically, the UC is generally formulated as a mixed-integer programming problem (MIP), which belongs to the class of NP-hard problems even if a single time period is considered \cite{anjos}. Moreover, the penetration of renewable generation---like wind and solar power---in power systems is rapidly deepening. Such growth in renewable generation, with their intermittent, variable and stochastic characteristics, puts increasing pressure on operators having to handle both asset failures and renewables' uncertainty and variability as part of the market-clearing process. Therefore, the development of computationally-efficient and robust methods to solve UC problems quickly and to optimality has been and continues to be a hot research topic \cite{YANG2021}.

Power system operators' experience and past research have shown that in UC problems whose solutions have to satisfy network constraints, only a small proportion of those  constraints can be potentially be binding---in other words, only a small proportion of lines can be congested \cite{Ardakani1, Ardakani2}. Considering only the potentially active transmission constraints significantly reduces UC solution times and can facilitate obtaining an optimal UC plan \cite{YANG2021}. The pruning process relies on identifying superfluous power flow constraints that can be safely eliminated from the original problem without jeopardizing its optimal solution nor its feasibility.\par

Various optimization-based approaches have been proposed for constraint screening for the UC problem. The iterative methods in \cite{Fu, Tejada-Arango, Chen2016} focus on solving a particular instance of the UC problem. Therefore, they can be practical in daily electricity market operations, but they encounter significant performance challenges in the context of short-term operations. Alternatively, the notion of {\textit{umbrella constraint}} was introduced when identifying redundant constraints which do not alter the feasibility region of the original UC problem when they are removed from the original problem \cite{Ardakani2}. Similarly, references \cite{Zhai, Roald, Zhang, Porras} use a \emph{bound tightening} technique \cite{Shchetinin} and solve two optimization problems for each transmission line in the power system over linear programming (LP) relaxations of the feasible region to remove as many redundant constraints as possible from the full UC formulation. Also, the authors of \cite{Weinhold} used \textit{Clarkson’s redundancy removal} algorithm that is based on LP and requires solving multiple maximization problems \cite{Szedlk2017RedundancyIL} for each possible line flow contingency constraint in a SCOPF. The drawback of \cite{Weinhold} is the ignorance of correlated patterns in historical samples of the load and renewable generation. Recently, an analytical approach was proposed by \cite{Ma} based on a constraint generation algorithm and heuristics to eliminate nonbinding constraints from UC problems on top of redundant constraints. \par

In addition, the variability and uncertainty of renewable power generation have introduced new challenges to the UC problem \cite{Mohandes}. The representation of uncertainty can take the form of scenarios or uncertainty sets for stochastic and robust optimization approaches, respectively. The computational burden of stochastic programming prompts to migrate towards tractable approaches for optimization under uncertainty, namely robust and chance-constrained optimization \cite{Wei,Venzke}. Moreover, renewable energy sources (RES) manifest cross-correlation over space and time. Predicting spatial and temporal scenarios has been of interest of power systems operation and planning community for a while now \cite{Golestaneh, Golestaneh2, Huo}. When modeling multivariate correlated scenarios, uncertainty sets take the form of boxes, polyhedral \cite{Golestaneh} and ellipsoidal sets \cite{Golestaneh2}. However, decision-making problems in power systems which include network constraints need different forms of uncertainty sets. For instance, integrating polyhedral uncertainty envelopes with a linear programming problem results in a linear programming problem, whereas the same problem with the ellipsoidal uncertainty sets is a second-order cone programming (SOCP) problem. Even though SOCP is convex, its nonlinearity is deemed a practical pitfall \cite{Golestaneh}.\par 

Within the context of constraint screening for the UC problems with deep penetration of correlated uncertainties, advancement in three aspects is required to close the gap between computationally tractable, robust, and optimal solutions to the UC constraint screening problem. First, it is necessary to design a computationally-tractable scalable optimization approach that can handle a high degree of uncertainty. Second, it is critical to define a robust and tractable uncertainty set that can capture stochastic dependence between different uncertainties---\emph{i.e.}, loads and renewable generation---, and be used as an input to a robust optimization approach \cite{Golestaneh}. Third, to yield optimal solutions in constraint screening problems means screening out not only redundant constraints with respect to the feasibility region but also as guided by the UC problem's objective function. In other words, a UC's non-active constraints can be part of its minimal feasible region, but they are never active constraints as they do not oppose the minimization the UC's objective function \cite{Pineda1,Porras}.  

From a methodological point of view, previous work \cite{Ardakani1, Ardakani2} developed mixed-integer linear programming (MILP) and LP formulations respectively, for umbrella constraint discovery (UCD) with a single residual demand parameter vector. The authors of \cite{Amir2} enhanced the original UCD notion from \cite{Ardakani1, Ardakani2} that has lighter computational demands and is better adapted when residual demand uncertainty is considered. The merit of the enhanced approach in \cite{Amir2} is the direct identification of umbrella constraints in lieu of the identification of non-umbrella constraints. 
References \cite{Roald, Aquino} considered uncertain residual demand parameters for robust optimization problems. However, these studies built box-shaped uncertainty sets for the univariate system net load (\emph{i.e.}, load less renewable generation) and ignored the inherent spatio-temporal couplings of renewable generation and demand. Although robust optimization techniques immunize the system according to very stringent potential events, it is easy to yield sub-economical outcomes \cite{Huo}. The conservativeness of a robust solution is directly related to the size of the uncertainty set \cite{Golestaneh}.\par

In light of this shortcoming, we extend the UCD formulation of Abiri-Jahromi and Bouffard \cite{Amir2}, by comprehensively capturing the spatial correlation of residual demands to identify the umbrella line flow limits in the UC problem. Consequently, we focus in this paper on multivariate polyhedral uncertainty regions. In this spirit, our prominent motivation is to construct a data-driven approach capable of generating computationally-tractable, conservative, and robust polyhedral uncertainty sets. We develop a model for the uncertainty set which we call \emph {data-driven polyhedral uncertainty set}. Compared to previous work, our proposed uncertainty envelope is less conservative than a conventional box-shaped uncertainty set \cite{Roald} and computationally cheaper than the convex hull of the uncertainty data \cite{Velloso}. Our proposal improves over a previous data-driven polyhedral uncertainty set approach \cite{Awadalla}. We demonstrate how it provides more conservative uncertainty coverage while, at the same time, improving significantly the computational performance of the constraint screening while guaranteeing the same outcome. 

Furthermore, we introduce a valid upper bound inequality constraint whose calculation is based on past UC solution runs as suggested by \cite{Porras} to eliminate further non-binding constraints. This allows the elimination of unrealistic generation schedules---with, for example, prohibitively high costs---as we run UCD. In contrast to the suggestion from \cite{Porras}, we include a conservativeness factor to mitigate model misspecification and infeasibility, a known weakness of \say{predict-then-optimize} approaches \cite{Elmachtoub}.

To this end, we tighten the LP relaxation of the UC by imposing the data-driven uncertainty set in the problem formulation first as an equality constraint, and second by a cost-driven inequality constraint. In doing so, we substantially increase the number of line flow constraints that should be eliminated from the original UC without jeopardizing its feasibility nor optimality. Furthermore, as the UCD algorithm can lend itself well to decomposition, we introduce a decomposition technique to further expedite its solution. 

The remainder of this paper is organized as follows. Section~II describes the formulation of the UC problem. Next in Section~IV, we revisit the notion of UCD. Section~IV describes how data-driven uncertainty sets and cost-driven constraints are obtained for tightening the UC relaxation. Section~V conducts case studies to show the effectiveness of our proposed approach. Finally, Section~VI presents our conclusions.

\section{Unit Commitment Model}
For expository purposes, we carry out our developments using a simplified single-period unit commitment \cite{Pineda1, Roald, Porras} according to the following simplifying assumptions:

\begin{itemize}
\item Single-period: UC is usually formulated as a multi-period problem that incorporates inter-temporal constraints typically associated to generation, for example minimum up and down times. However, since this work focuses on investigating the impact of residual demand spatial correlation on network reduction, we prefer to investigate this solely by considering a single-period UC  \cite{Pineda1, Roald}. 
\item DC power flow: The power flows in the transmission  lines are estimated via a dc approximation by using power transfer distribution factors (PTDF) to keep the model linear. The PTDF of line $l$ with respect to node $n$ is
denoted as $h_{l n}$. Besides, $f_{l}^{\max}$ represents the maximum flow capacity of transmission line $l$. The number of buses and lines are denoted by $N$ and $L$, respectively. 
\item Generation portfolio: Each generating unit $m$ is characterized by a minimum and a maximum power output which are denoted as $g_{m}^{\min}$ and $g_{m}^{\max}$, respectively. 
\item Residual demand: The residual demand at bus $n$, $d_n$, is a net load (demand less non-dispatchable renewable generation). Without loss of generality, we assume that residual demands are multivariate normally distributed correlated random variables.     
\item No contingencies: We assume that all generators and lines are fully operational and therefore security constraints are neglected.
\end{itemize}

\subsection{Problem Formulation}
The optimization problem corresponding to this simplified UC is a MILP problem and formulated as in \cite{Pineda1}:
\begin{equation}
\min _{u_{m}, g_{m}, q_{n}} \sum_{m \in \mathcal{M}} c_{m} g_{m}
\label{eq:P1a}
\end{equation}
Subject to:
\begin{align}
q_{n} &= \sum_{m \in \mathcal{M}_{n}} g_{m}-d_{n},  & \forall n \in\mathcal{N} \label{eq:P1b} \\
\sum_{n=1}^N q_{n} &= 0 & \label{eq:P1c}\\
u_{m} g_{m}^{\min} & \leq g_{m} \leq u_{g} g_{m}^{\max}, & \forall m \in\mathcal{M} \label{eq:P1d} \\
-f_{l}^{\max } & \leq \sum_{n=1}^N h_{l n} q_{n} \leq f_{l}^{\max }, & \forall l \in\mathcal{L} \label{eq:P1e} \\
u_{m} &\in\{0,1\}, & \forall m \in\mathcal{M} \label{eq:P1f}
\end{align}
Decision variables include the commitment status of the generating units $u_{m}$, the power output schedules  $g_{m}$, the net power injections at each node $q_{n}$. The objective function \eqref{eq:P1a} minimizes the total generation cost. Constraint \eqref{eq:P1b} computes the net injected power at each node, while constraint \eqref{eq:P1c} ensures power balance in the system. Constraints \eqref{eq:P1d} and \eqref{eq:P1e} respectively enforce limits on generator outputs and power flows on transmission lines using PTDFs. Finally, \eqref{eq:P1f} requires that the on/off status of generators are binary ($u_m = 0$ if off, $u_m = 1$ if on).

\section{Umbrella Constraint Discovery}
The \emph{umbrella constraints} \cite{Ardakani1} of an optimization problem are the constraints that, if removed, change the feasibility region of that problem. The \emph{umbrella set} of an optimization problem is the set of constraints containing the fewest constraints preserving the original optimization problem feasibility region.
\subsection{Identification of Umbrella Network Constraints in UC}
The UC problem \eqref{eq:P1a}--\eqref{eq:P1f} can be made significantly easier to solve if constraints \eqref{eq:P1e} with no impact on the optimal UC plan are removed \cite{Ardakani1}. In this paper, we tailor the UCD problem with the objective of favoring the identification of potentially active transmission lines power flow constraints in the unit commitment problem rather than identifying the complete umbrella set.
Here, in its original incarnation UCD, see for example \cite{Ardakani1}, can identify the set of umbrella line flow limits constraints. That minimal set of constraints is necessary and sufficient to characterize the feasible region of the original UC problem. UCD is an iterative algorithm which, at each iteration, solves the optimization problem \eqref{eq:P2a}--\eqref{eq:P2l} stated next. Each iteration finds the set of line constraints forming one of the vertices of the feasible region of the original UC problem. Once all vertices have been found, the algorithm terminates. Specifically at each iteration, we solve for the binary vectors $v_{l}^{\pm} \in\{0,1\}^{L}$, and continuous vectors $g \in \mathbb{R}^{M}$, $q \in \mathbb{R}^{N}$, $d \in \mathbb{R}^{N}$ and $z^{\pm} \in \mathbb{R}_{+}^{L}$.
\begin{align}
\min & \sum_{l \in \mathcal{L}} (v_{l}^{+} + v_{l}^{-}) & \label{eq:P2a} \\
\intertext{Subject to:}
q_{n} &= \sum_{m \in \mathcal{M}_{n}} g_{m}-d_{n},  & \forall n \in\mathcal{N} \label{eq:P2b} \\
\sum_{n=1}^N q_{n} &= 0 & \label{eq:P2c}\\
u_{m} g_{m}^{\min} & \leq g_{m} \leq u_{g} g_{m}^{\max}, & \forall m \in\mathcal{M} \label{eq:P2d} \\
\sum_{n=1}^N h_{l n} q_{n} & \leq f_{l}^{\max }, & \forall l \in\mathcal{L} \label{eq:P2ee} \\
-\sum_{n=1}^N h_{l n} q_{n} & \leq f_{l}^{\max }, & \forall l \in\mathcal{L} \label{eq:P2e} \\
\sum_{n=1}^N h_{l n} q_{n} + z_l^{+} & \geq f_{l}^{\max }, & \forall l \in\mathcal{L} \label{eq:P2f} \\
-\sum_{n=1}^N h_{l n} q_{n} + z_l^{-} & \geq f_{l}^{\max }, & \forall l \in\mathcal{L} \label{eq:P2ff} \\
v_{l}^{+} - \frac{z_{l}^{+}}{\Omega} & \geq 0, & \forall l \in\mathcal{L} \label{eq:P2g} \\
v_{l}^{-} - \frac{z_{l}^{-}}{\Omega} & \geq 0, & \forall l \in\mathcal{L} \label{eq:P2gg} \\
z_{l}^{+}, z_{l}^{-} &\geq 0, & \forall l \in\mathcal{L} \label{eq:P2h}\\
v_{l}^{+}, v_{l}^{-} &\in\{0,1\}, & \forall l \in\mathcal{L} \label{eq:P2i}\\
0 &\leq u_{m} \leq 1, & \forall m \in\mathcal{M} \label{eq:P2k}\\
d_{n}^{\min } & \leq d_{n} \leq d_{n}^{\max }, & \forall n \in\mathcal{N} \label{eq:P2l} 
\end{align}
where $\Omega$ is a large positive number. Here, the binary variables $v^{\pm}_l$ take the value of 0 if one of the flow limits associated with line $l$ are umbrella ($v_l^+ = 0$ if the upper flow limit is umbrella or $v_l^- = 0$ if the lower flow limit is umbrella). Otherwise, $v^{\pm}_l$ are set to 1.

The objective function \eqref{eq:P2a} aims to minimize the sum of the binary variables, $v^\pm_{l}$, by finding the maximum number of line flow constraints that can hit or surpass their capacity.

The set of constraints \eqref{eq:P2b}--\eqref{eq:P2d} from the UC, controls the decision variables $q_{n}$ and $g_{m}$. Each line flow constraint in the blocks of constraints \eqref{eq:P2ee} and \eqref{eq:P2e} is paired with a constraint in the block \eqref{eq:P2f} and \eqref{eq:P2ff}. By inspection, the auxiliary variables $z_{l}^{\pm}$  can be equal to zero only if there is $\sum_n q_n$ satisfying both \eqref{eq:P2ee} and \eqref{eq:P2f} for the upper flow limits, or \eqref{eq:P2e} and \eqref{eq:P2ff} for lower flow limits. As a result, the binary variable $v_{l}^{\pm}$ = 0. Conversely, $z_{l}^{\pm}$ has to be positive and the binary variables $v_{l}^{\pm}$ = 1 as required by \eqref{eq:P2h}. Furthermore, the constraints for which $v_{l}^{\pm}$ = 1 have to intersect at the same value of $\sum_n q_n$. Since $\sum_n q_n$ is an intersection of line constraints over LP-relaxation of the feasible set of the UC problem, it is therefore a vertex of this set. Additionally, the vector of residual demand at each bus $n$ is turned into a vector of decision variables as in \eqref{eq:P2l}.

Following the first iteration, which revealed the vertex with the most intersecting constraints, the next step is to pinpoint the other vertices that have the same or fewer intersecting umbrella constraints. This is done while the previously discovered umbrella constraints are removed from the search by setting their respective binary variable $v_l^\pm$ equal to 1. We terminate the search when there are no more umbrella constraints to identify. The UCD procedure is summarized in Algorithm~\ref{alg:one}.
 
The drawback of the UCD algorithm relates to the modeling of residual demand range limits in constraint \eqref{eq:P2l}. This result will be sub-optimal because it ignores the spatial correlation among the residual demands. Therefore, we propose next a data-driven UCD (D-UCD) where the residual demand vector is characterized using polyhedral uncertainty sets as described in Section III, rather than a simple box in $N$-dimensional space. We will refer to the UCD framework described in \eqref{eq:P2a}--\eqref{eq:P2l} as the base UCD (B-UCD) approach. In addition, note that UCD solution benefits from decomposition; the interested reader can view its details in Appendix~A.

\begin{algorithm}
\caption{Umbrella Constraint Discovery}\label{alg:one}
\KwData{Network and generators data, historical net-load}
\KwResult{Non redundant constraints.}
\While{$\sum_{l} v_l^{+} + v_l^{-} \neq 2L $}{
 Solve UCD \eqref{eq:P2a}--\eqref{eq:P2l}\;
  \eIf{($v_l^{+}=0$ OR $v_l^{-}=0$)}{
   Set $v_l^{+}=1$ OR $v_l^{-}=1$\;
   go  to step 2\;
   }{
   go  to step 9\;
  }
  }
\end{algorithm}


\section{Tightening Umbrella Constraint Discovery}
As argued in the Introduction, our goal is to eliminate both redundant and inactive constraints from the UC problem formulation with aim to run a computationally-lighter unit commitment. At the same time, we have to target a formulation is able to map adequately net load uncertainty with \emph{as few constraints as possible}. As a first step to accomplish this, we compute uncertainty sets of residual demands as a function of historical records of net load and their initial forecasts.  

\subsection{Data-Driven Polyhedral Uncertainty Sets}
Without loss of generality, we focus on capturing the spatially correlated uncertain residual demand. Inspired by \cite{Pachamanova, Golestaneh}, we develop two formulations of data-driven polyhedral uncertainty sets by leveraging principal component analysis (PCA) \cite{Huo, Berizzi}. PCA is applied to historical time series of net load. All the input time series have the same length of $T$, and they are synchronized and evenly spaced in time (\emph{e.g.} one hour intervals).  \par

We denote a matrix ${W} \in {\mathbb{R}}^{T \times N}$ whose elements $w_{n t}$ are the time series of \emph{observed} residual demand at bus $n \in \mathcal{N} = \{1,...,N\}$ for each time instance $t \in \mathcal{T} = \{1,...,T\}$. Similarly, we define the matrix $\mu \in {\mathbb{R}}^{T \times N}$ whose elements $\mu_{n t}$ are the time series of past \emph{forecasted} net demand at bus $n \in \mathcal{N}$ for each past time instance $t \in \mathcal{T}$.

Using $\mu$, we obtain the centered data matrix $W_{c}$ \cite{Huo, Berizzi}, whose contents are the residual demand forecast errors at all nodes $n \in \mathcal{N}$ and times $t \in \mathcal{T}$
\begin{equation}
W_{c}=W-\mu
\end{equation}
Assuming residual demand forecast errors are unbiased\footnote{In the case where errors are unbiased, one would need to calculate the biases at each node and then remove them from $W_c$.}, its spatial forecast error covariance matrix $\Sigma \in {\mathbb{R}}^{N \times N}$ is approximated by
\begin{equation}
\Sigma=\frac{1}{T-1} {W}_{c}^\top {W}_{c}
\end{equation}

PCA is performed by conducting an eigenvalue decomposition of the covariance matrix. We let the columns of an $N \times N$ matrix $V$ and the diagonal entries of another $N \times N$ matrix $\Lambda$ represent, respectively, the orthonormal eigenvectors and the eigenvalues of $\Sigma$. Here, the diagonal elements of $\Lambda$ are ordered such that $\lambda_{11} \geq \lambda_{22} \geq \cdots \geq \lambda_{NN}$, while the columns of $V$ are arranged such that its $n$th column (eigenvector) is associated with the $n$th eigenvalue, $\lambda_{nn}$.

Thus, we can project the data contained in $W_c$ onto the eigenvectors $k = 1,\ldots,N$ of the covariance matrix \cite{Berizzi}
\begin{equation}
{Z_k}= {W}_{c} V_k
\end{equation}
where $Z_k$ is a $T \times 1$ vector of data which has been projected onto $V_k$, the $k$th principal component of $\Sigma$.

Next, let us find the coordinates of the extrema of each data projection $Z_k$
\begin{equation}
    \bar{\mathcal{S}}_k = \arg\max_{t \in \mathcal{T}} \| Z_k \|^2
\end{equation}
that is $\bar{\mathcal{S}}_k$ is the data point projected along principal component $V_k$ which is the furthest away from the origin. Keeping a conservative approach, we will assume data can range between $-\bar{\mathcal{S}}_k$ and $\bar{\mathcal{S}}_k$ along the principal component $V_k$. Moreover, we can argue that the extrema of the original data points can be reconstructed using the $K$ principal components and their data projections $Z_k$ recentered on the net load forecast $d^0 \in \mathbb{R}^N$
\begin{align}
    \mathcal{S}^+_k &= d^0 + \bar{\mathcal{S}}_k \\
    \mathcal{S}^-_k &= d^0 - \bar{\mathcal{S}}_k
\end{align}

Considering that typically only the first few dominant principal components are sufficient to describe accurately the original data's uncertainty, it is common practice to limit the number of principal components to $K < N$. A first \emph{data-driven polyhedral uncertainty set} (DPUS) as proposed in \cite{Pinar}, is
\begin{align}
P_1(\mathcal{S},K) =& \Bigg \{\mathbf{E} \in \mathbb{R}^N \mid \mathbf{E} = 
\sum_{k=1}^{K} \left ( \omega_{k} \mathcal{S}^+_{k} + (1-\omega_{k}) \mathcal{S}^-_{k} \right ) \nonumber \\
& \qquad  0 \leq \omega_{k} \leq 1, k \in \{1, \ldots, K\} \Bigg \} \label{eq:Pus1} 
\end{align}

\begin{figure}
\begin{center}
\includegraphics[width=8.0cm]{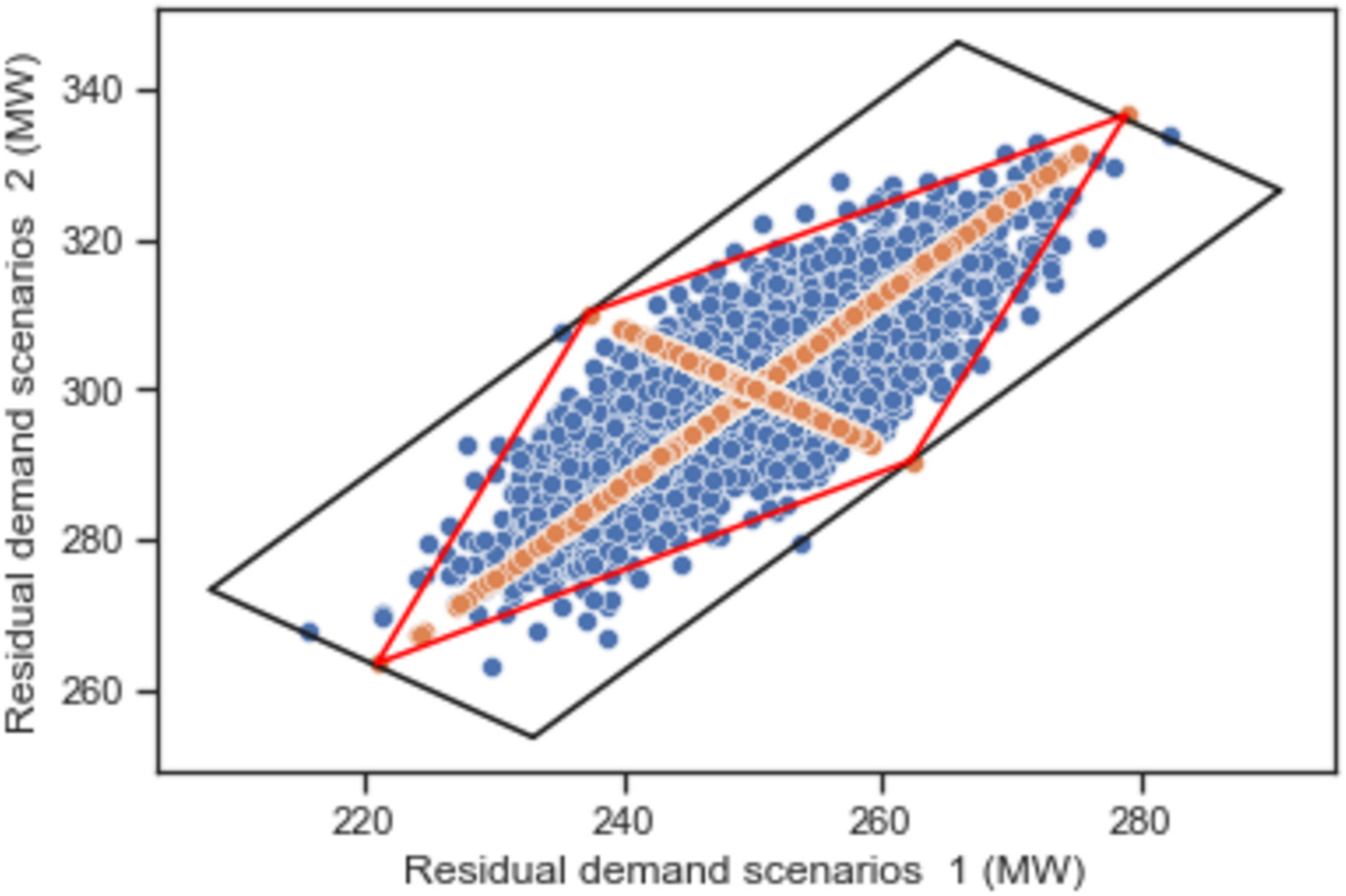}
\caption{Schematic illustration of polyhedral uncertainty envelope with $d^0 = (250, 300)$~MW and $K=2$.}
\label{Fig_Usets}
\end{center}
\end{figure}

Fig.~\ref{Fig_Usets} illustrates the construction of a polyhedral uncertainty set in a two dimensional uncertainty space. Fig.~\ref{Fig_Usets} shows $P_1(\mathcal{S},K)$ for $K=2$ in black. The orange dots represent the original data as it is projected along the two principal components of historical forecast error data. By inspection, we see that this polyhedral set encloses all original data points.

In line with previous proposals \cite{Awadalla}, we can define an alternative polyhedral uncertainty set
\begin{align}
P_2(\mathcal{S},K) =& \left\{\mathbf{E} \in \mathbb{R}^N \mid \mathbf{E} = \sum_{k=1}^{K} \left ( \omega^+_{k} \mathcal{S}^+_{k} + \omega^-_{k} \mathcal{S}^-_{k} \right ), \right . \nonumber \\
& \qquad \sum_{k=1}^{K} \left ( \omega^+_{k} + \omega^-_{k} \right ) = 1, \nonumber \\
& \qquad  0 \leq \omega^+_{k} \leq 1, \; 0 \leq \omega^-_{k} \leq 1, \nonumber \\
& \qquad k \in \{1,\ldots,K\} \Bigg \}\label{eq:Pus2} 
\end{align}
where we notice that $P_2(\mathcal{S}, K) \subseteq P_1(\mathcal{S}, K)$.

The set $P_2(\mathcal{S},K)$ represents the smallest convex set that contains every data point projected onto the $K$ retained principal components. Moreover, that DPUS is a convex hull of the extrema of the retained $Z_k$ data projections, where we define $\mathcal{S} = \cup_k ( \mathcal{S}_k^+ \cup \mathcal{S}_k^-)$. In Fig.~\ref{Fig_Usets}, we see the historical data, represented by blue dots, and the data projected onto two of its principal components (orange dots). By inspection, the rhombus-shaped red envelope, whose principal axes correspond to the principal components of the data, encapsulates the vast majority of the original data. 

Later in Section~V, we will examine the pros and cons of these two net load uncertainty representations as applied to the UC problem and the reduction of its number of constraints.

\subsection{Economics-Driven UCD (ED-UCD)}
Another relevant limitation of conventional UCD as seen in \eqref{eq:P2a}--\eqref{eq:P2l} is that it solely seeks to establish the minimal set of constraints required to describe the feasible space of the UC. It is incapable of providing information regarding which of the umbrella constraints could become active as we solve the UC. We argue that having such information ahead of UC solution would be a valuable asset in UC solution time reduction.

\begin{figure}
\begin{center}
\includegraphics[width=7.0cm]{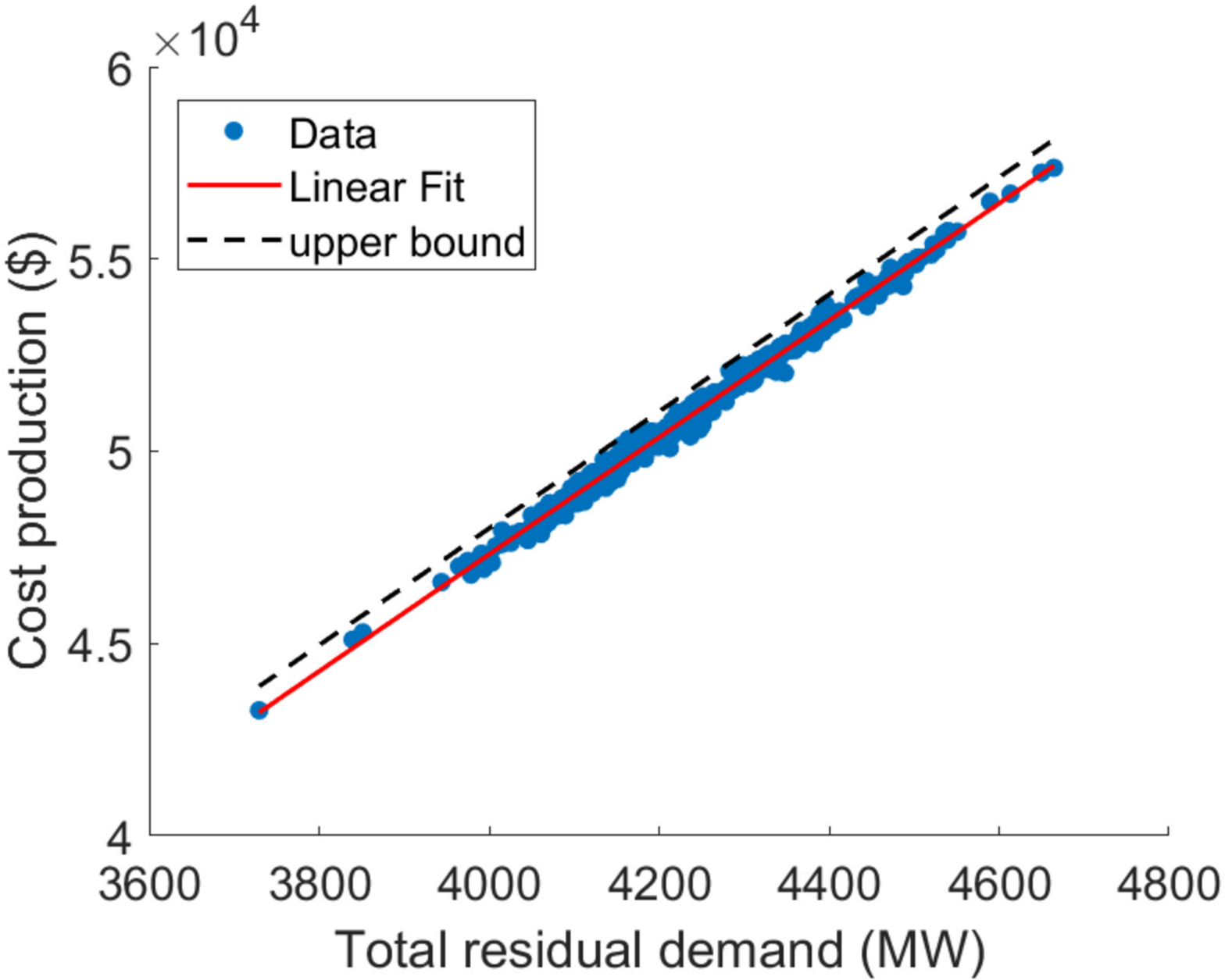}
\caption{Observed production cost against system net load.}
\label{Fig_Usets2}
\end{center}
\end{figure}

To illustrate this, consider Fig.~\ref{Fig_Usets2} that shows how the total production cost of a UC problem formulation as a function of the total residual demand for 300 solution instances of a given UC problem. Clearly, the positive correlation between the operating cost and total net demand drives towards modeling this relationship using a linear regression technique.\footnote{We note that a family of such curves may be necessary to cover a wider range of net demands. One would expect that these taken together would form a piecewise linear convex function of the total net demand.} We argue that all UC outcomes in a given range of net demand would have their corresponding production costs bounded above the dashed line shown on Fig.~\ref{Fig_Usets2}.

Therefore, we propose that if we are to add new constraints to the UC, which we cast as follows \cite{Porras}:
\begin{equation}
 \sum_{m \in \mathcal{M}} c_{m} g_{m} \leq (1 + \Delta \sigma) a_{0} + (1 + \Gamma) b_{0} D \label{eq:ED1}
\end{equation}
\begin{equation}
D=\sum_{n=1}^N d_n
\label{eq:ED2}
\end{equation}
\begin{equation}
D^{\min }  \leq D \leq D^{\max }
\label{eq:ED3}
\end{equation}
and perform UCD on the resulting problem, we would be able to identify which of the original problem's umbrella constraints are most likely to be active at the UC's own optimum. The premise here is that \eqref{eq:ED1} should be intersecting the active umbrella constraints of the original UC. Obviously, here there is some tuning to be carried out in determining the parameters $a_{0}$ and $b_{0}$ of the proposed cost upper bound, while $\Delta \sigma$ and $\Gamma$ are user-specified conservativeness parameters.

The upper bound uses a basic linear fitted model $a_{0}+b_{0}D$. The minimum and maximum aggregate net load for each is denoted by $D^{\min}$ and $D^{\max}$, respectively. The upper bound can impact the screening outcomes in terms of expected eliminated constraints.

Through the selection of the value of $\Delta \geq 0$, which multiplies the standard deviation $\sigma$ of its underlying data with respect to the best fit line, it is possible to push up on the cost upper bound. This way it is possible to capture most if not all prior cost observations. For example, with $3\sigma$ one will typically capture production costs of all prior observed instances as shown by the dashed line in Fig.~\ref{Fig_Usets2}. On the other hand, if $\Delta$ is set too high, there will be a risk that the cost upper bound \eqref{eq:ED1} is in fact found to be redundant when running UCD on the augmented UC constraint set and, thus, be of little value. Moreover, we note that if $\Delta$ is too small, we run the risk that \eqref{eq:ED1} renders the augmented UCD infeasible.

Second, we add the factor $\Gamma \geq 0$ to the linear model to further trade-off between the number of retained umbrella constraints and the risk of excluding an umbrella constraint which may become active in the UC. The consequence of excluding potentially active constraints will be that UC solutions may be infeasible since one or more UC constraints are not satisfied by virtue of having been excluded by the augmented UCD. In fact, by setting $\Gamma > 0$ one is lowering the risk of the augmented UCD weeding out constraints that need to be considered in the UC.


Finally, we note that the economic-driven upper bound can be extended using a piecewise linear set of constraints to capture residual demand and cost data over different ranges of net demand as suggested by \cite{Porras}.


Next, we will illustrate how the combination of data-driven polyhedral uncertainty sets and ED-UCD can reduce dramatically the computational effort required to solve robust UC problems. Table~\ref{table1} summarizes the various tightened UCD problems which we will be comparing in the following section.


\begin{table}
\captionsetup{font=footnotesize}
\caption{\sc{Constraints Screening Methods}} 
\centering 
\footnotesize{
\begin{tabular}{c c} 
\hline
Method          & Screening optimization problem   \\ [0.5ex] 
\hline 
B-UCD           & \eqref{eq:P2a} s.t \eqref{eq:P2b}--\eqref{eq:P2l}      \\ [0.5ex]
D1-UCD          & \eqref{eq:P2a} s.t \eqref{eq:P2b}--\eqref{eq:P2k}, \eqref{eq:Pus1}      \\ [0.5ex]
D2-UCD          & \eqref{eq:P2a} s.t \eqref{eq:P2b}--\eqref{eq:P2k}, \eqref{eq:Pus2}    \\ [0.5ex] 
ED-UCD          & \eqref{eq:P2a} s.t \eqref{eq:P2b}--\eqref{eq:P2l}, \eqref{eq:ED1}--\eqref{eq:ED3}  \\ [0.5ex]  
ED+D1-UCD       & \eqref{eq:P2a} s.t \eqref{eq:P2b}--\eqref{eq:P2k}, \eqref{eq:Pus1}, \eqref{eq:ED1}--\eqref{eq:ED3}  \\ [0.5ex] 
\hline 
\end{tabular}}
\label{table1} 
\end{table}

\section{Case Study}
\subsection{Benchmark Approach (BA): Roald's Method}
This method was proposed in \cite{Roald}; it is based upon the solution of one maximization and one minimization for each transmission line $\hat{l}$, these two optimizations are jointly formulated as 

\begin{equation}
\max / \min \sum_{n=1}^N h_{\hat{l} n} q_{n}
\label{eq:P3a}
\end{equation}
Subject to:
\begin{align}
\eqref{eq:P2b} & -\eqref{eq:P2e}, \; \eqref{eq:P2k}-\eqref{eq:P2l} & \label{eq:P3b}
\end{align}
In short, problem \eqref{eq:P3a}--\eqref{eq:P3b} seeks to maximize/minimize the power flow through each transmission line $\hat{l}$ over an LP-relaxation of the feasible region of the UC problem. If the maximum (minimum) limit for the flow of line $\hat{l}$ given by the objective function does not reach line capacity limit $f_{l}^{\max }$, then the upper (lower) line constraint is flagged as redundant.

\subsection{Procedure for Constructing Correlated Residual Demand Time Series}
We generate $N$ synthetic spatially-correlated residual demand time series of length $T$, which are then consigned to matrix $W$. They consist of historic residual demand forecasts $\mu \in \mathbb{R}^{T \times N}$, which correspond to the nominal demand values from the data sets in \cite{data2}. These are superimposed with zero-mean normally-distributed forecast errors with spatial correlation given by a covariance matrix $\Sigma$. Here, we take the approach outlined in \cite{Pena-Ordieres}, where errors are assumed to be proportional to forecasts and whose variance and correlation are adjusted with an uncertainty level parameter. Thus, this parameter controls the magnitude of net load forecast errors. We utilize the exact data generation approach proposed by \cite{Pena-Ordieres} which includes a random process in modeling residual demands' correlation. First, we generate a positive definite matrix $C=\widehat{C} \widehat{C}^{\top}$ where each element of the matrix $\widehat{C}$ is a sample randomly drawn from a uniform distribution with support in $[0,1]$. Then, to obtain a positive definite covariance matrix in which the diagonal elements are $c_{nn} = (\eta d^0_n)^2$, and off-diagonals
\begin{equation}
\sigma_{n m}=\eta^2 \frac{c_{n m}}{\sqrt{c_{n n} c_{mm}}} d^0_n d^0_{m}, \quad \forall n, m \in \mathcal{N}, n \neq m
\end{equation}
where $c_{nm}$ are the $nm$th elements of the matrix $C$, and $d_n^0$ is the net demand forecast at node $n$. Finally, we generate $T$ = 8760 nodal net load vectors using the described approach. Of the 8640 values of residual demand generated, 7200 are randomly selected to obtain the data-driven cost upper bound \eqref{eq:ED1} and uncertainty sets \eqref{eq:Pus1} and \eqref{eq:Pus2}, and the remaining 1440 instances are used for testing and investigating UC performance. \par


\subsection{Performance Evaluation}
The procedure to measure the performance of the method described in Section~III and its coupling with Section~IV to remove redundant line flow constraints---with the overarching objective of reducing the UC solution computational effort---is run as suggested by \cite{Pineda1}:

\begin{enumerate}
\item Given the historical data, set up net load uncertainty representations (as described in Section~III). Then, determine the transmission constraints that can be eliminated according to each approach described in Table.~\ref{table1}. 
 \item Record the computational time needed to screen out the redundant network constraints using each screening approach. We consider that the computational time of each approach is given as the sum of the time required to run each iteration for different versions of the UCD algorithm. On the other hand, the benchmark approach runs in a sequential manner for transmission line constraint screening. 
 \item Solve the reduced UC problem on the set of unseen time periods without the superfluous constraints identified in Step 1. 
 \item Use the binary commitment variables obtained in Step 3) as a warm start solution and solve the unit commitment problem including all constraints.
 \item Assess the performance of the screening method in terms of (i) the percentage of retained network constraints from each screening approach in Step 1, (ii) the computational time required to solve each screening approach, (iii) the computational time required to solve the reduced UC problem in Step 3) with respect to the full UC formulation, and (iv) whether or not all necessary and sufficient constraints needed to solve the UC have been retained. The latter performance measure applies to ED-UCD and ED+D1-UCD since the application of the production cost upper bound may discard potentially active UC constraints out of the original umbrella set, as explained in Section~IV.-B.
\end{enumerate}

\begin{table}[t!]
\captionsetup{font=footnotesize}
\caption{\sc{Description of test power systems}} 
\centering 
\footnotesize{
\begin{tabular}{c c c c c} 
\hline
System & \# Nodes & \# Generators & \# Lines  \\ [1.0ex] 
\hline 
IEEE-RTS-73             & 73   & 96 & 120 \\ [1.0ex]
IEEE-118                & 118  & 19 & 186 \\ [1.0ex] 
\text{CASE500\_pserc}   & 500  & 49 & 733 \\ [1.0ex]
\hline 
\end{tabular}}
\label{table2}
\end{table}

\subsection{Experimental Setup}
The B-UCD, D1-UCD and D2-UCD screening approaches and the UC are formulated as MILP problems. We test our algorithm in two standard IEEE test networks, namely IEEE-RTS-73 and IEEE-118 test systems \cite{data1}. Also, another test case is adopted from the IEEE PES PGLib-OPF v17.08 benchmark library \cite{data2} which is called \text{CASE500\_pserc}. All the technical data related to these systems are available in \cite{data1}, and their main features are listed in Table.~\ref{table2}. For these medium size test networks, the solution optimality gap was set to 0\%. For the benchmark approach, we run $2L$ optimizations in a sequential manner. The calculations have all been performed using GAMS and the CPLEX MILP solver. The computer used is equipped with an Intel Core i7 3.10~GHz processor and 16~GB of RAM. \par

\begin{table}
\captionsetup{font=footnotesize}
\caption{\sc{Solution time for redundant constraint screening for medium size networks}} 
\centering 
\footnotesize{
\scalebox{0.9}{
\begin{tabular}{c c c c c} 
\hline
& \multicolumn{4}{c}{Screening time (s)} \\
\hline 
Method  & BA & B-UCD & D1-UCD & D2-UCD   \\ [0.5ex] 
     &  & (\% change) & (\% change) & (\% change) \\[0.5ex]
\hline 
IEEE-RTS-73  &    162.5       & 36.7           & 45.3            & 58.1                \\[1ex]
         & -               & ($-$77.4\%)      & ($-$72.1\%)       & ($-$64.3\%)              \\[1ex]
IEEE-118   & 253.5             & 34.3           & 46.4           & 48.6                  \\[1ex] 
         & -               & ($-$86.4\%)      & (-81.7\%)      & ($-$80.8\%)                \\[1ex]                       
\hline 
\end{tabular}}}
\label{table3} 
\end{table}

\begin{figure}
\centering
\begin{subfigure}{0.232\textwidth}
    \includegraphics[width=\textwidth]{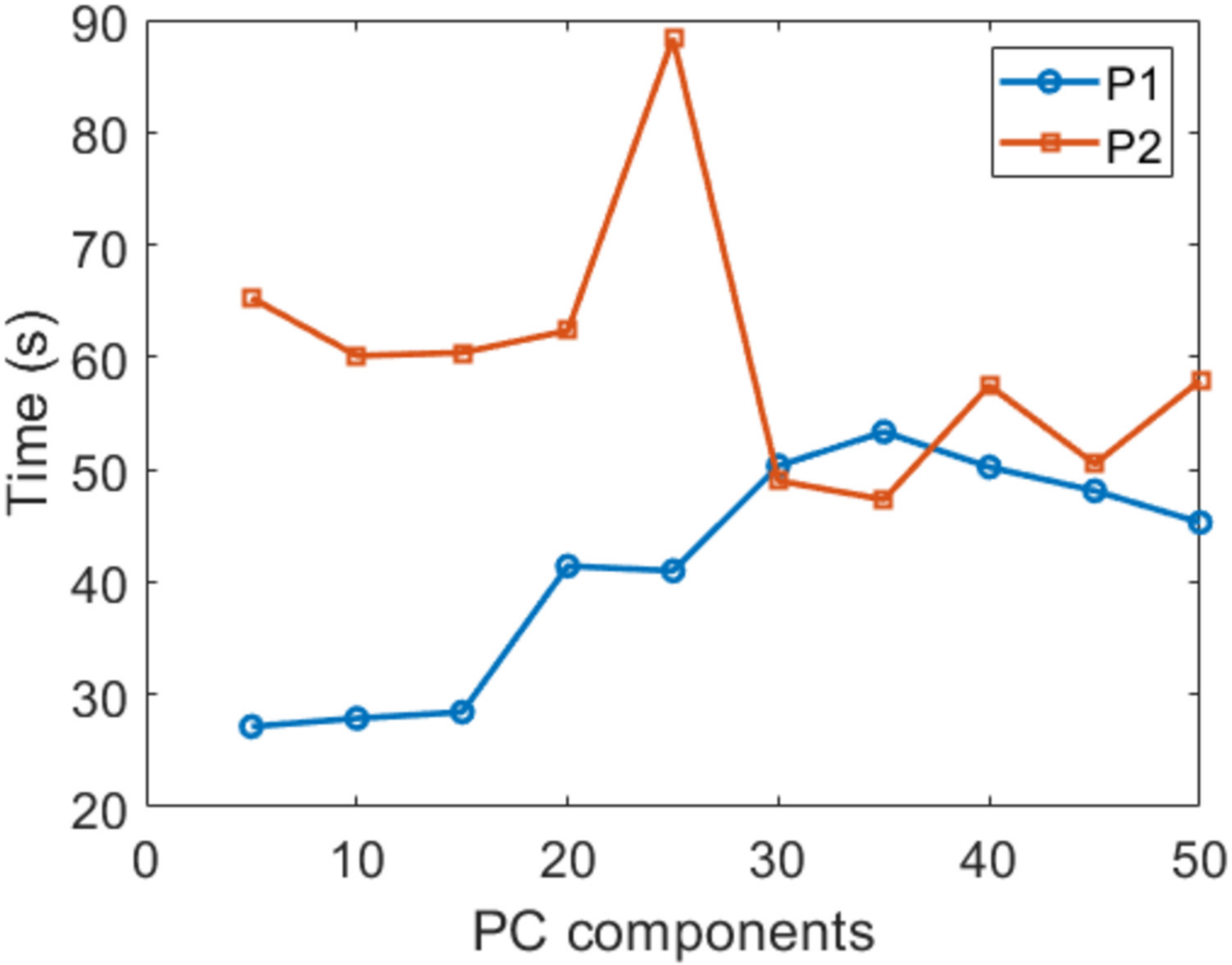}
    \caption{}
    \label{fig:04a}
\end{subfigure}
\hspace{0.5em}
\begin{subfigure}{0.232\textwidth}
    \includegraphics[width=\textwidth]{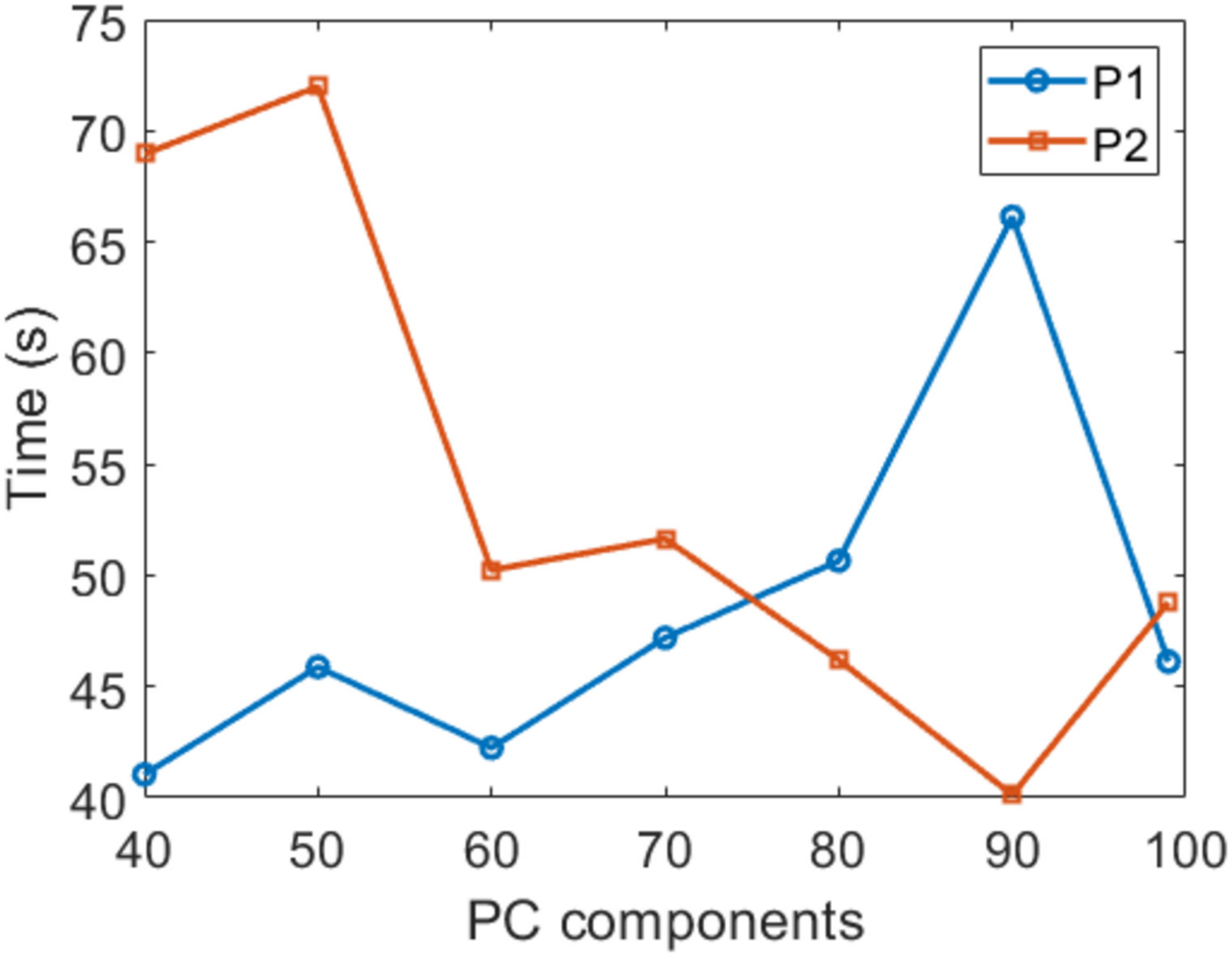}
    \caption{}
    \label{fig:04b}
\end{subfigure}
\begin{subfigure}{0.232\textwidth}
    \includegraphics[width=\textwidth]{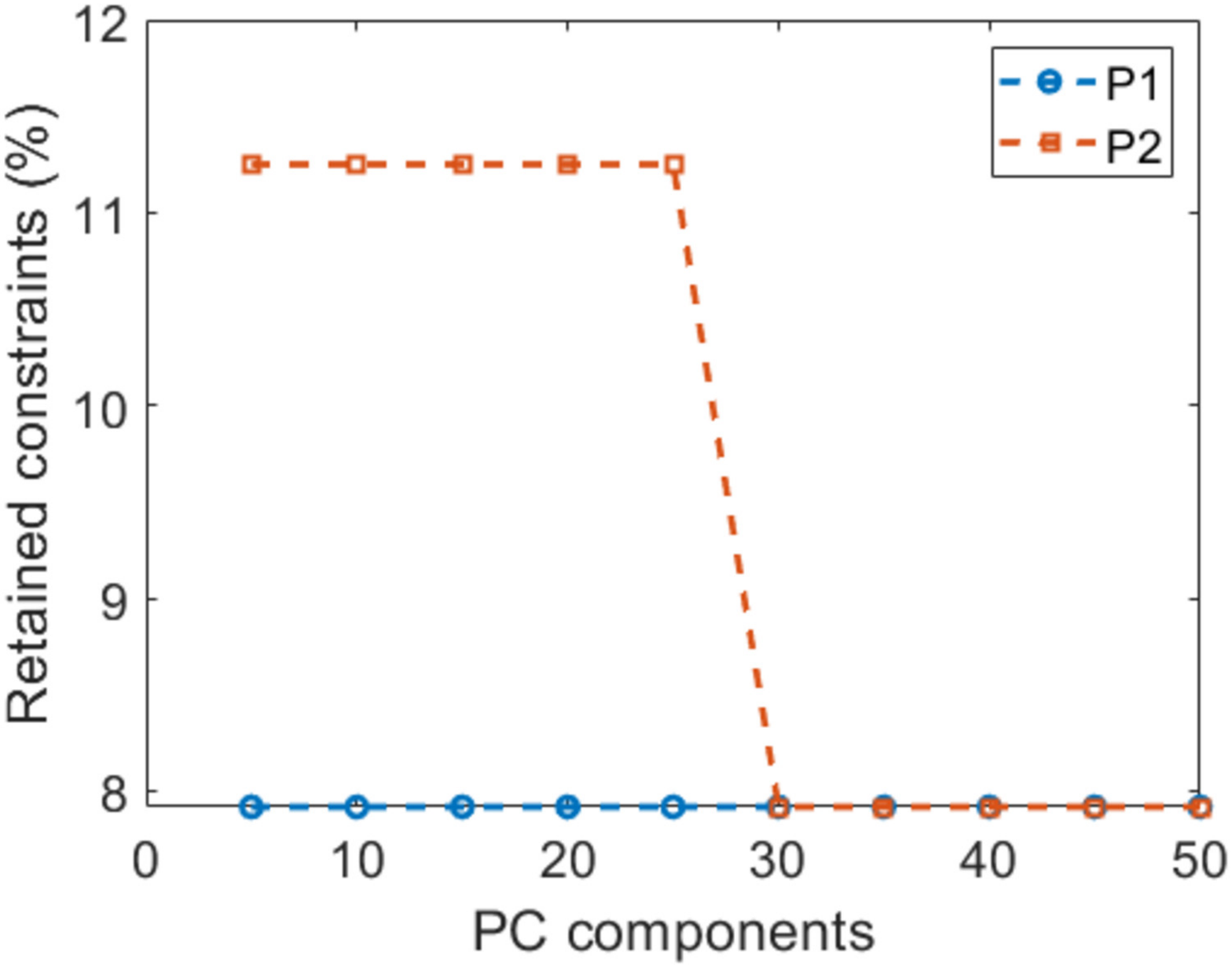}
    \caption{}
    \label{fig:04c}
\end{subfigure}
\hspace{0.5em}
\begin{subfigure}{0.232\textwidth}
    \includegraphics[width=\textwidth]{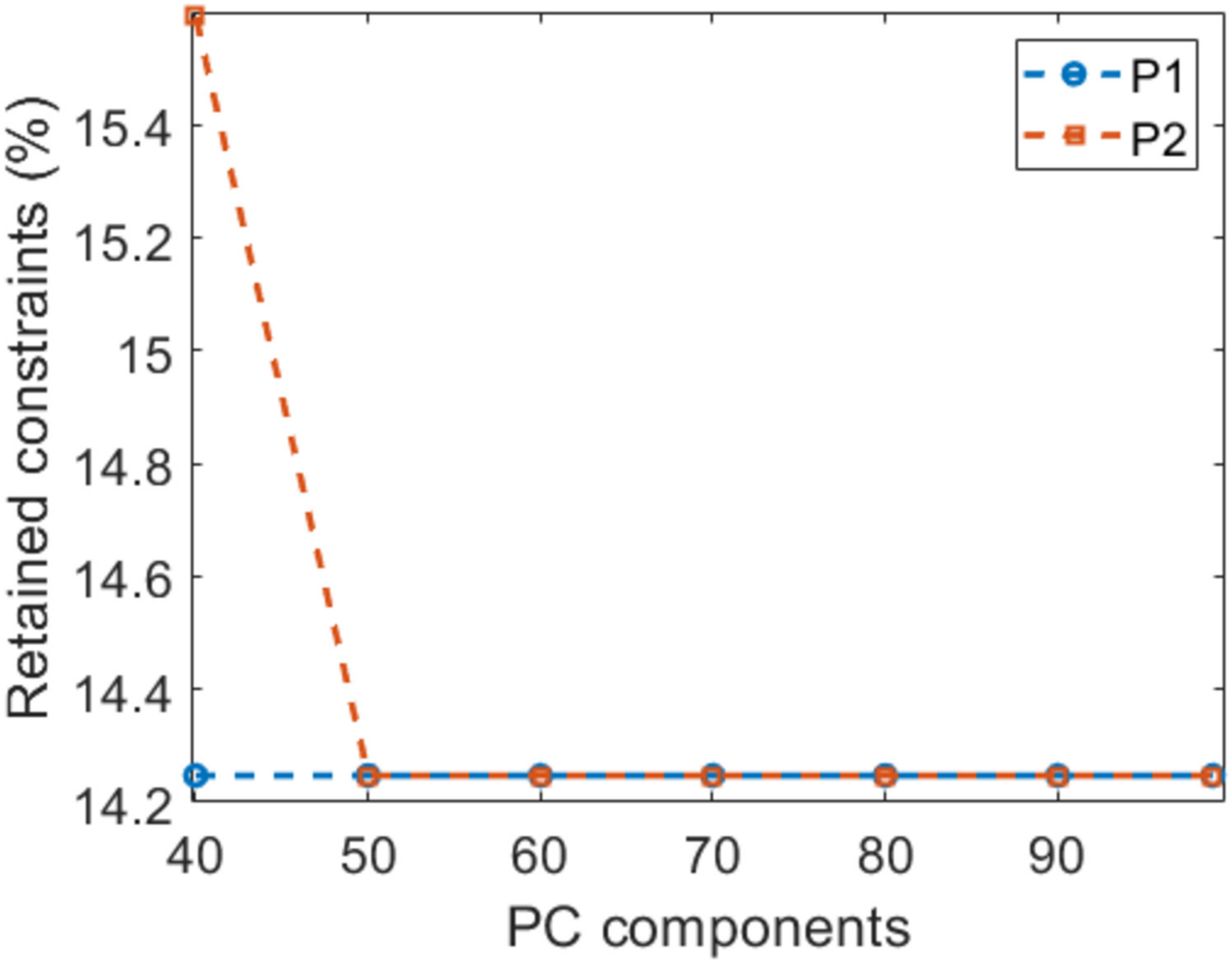}
    \caption{}
    \label{fig:04d}
\end{subfigure}
\hfill
\caption{(a) The impact of $K$ on screening time in IEEE-RTS-73; (b) The impact of $K$ on screening time in IEEE-118. (c) The impact of $K$ on retained constraints in IEEE-RTS-73; (d) The impact of $K$ on retained constraints in IEEE-118.}
\label{fig:04}
\end{figure}

\subsection{Redundant Constraints' Screening Results}
\subsubsection{Medium-Size Networks}
In this subsection, we provide simulation results for the two medium-size test systems provided in Table~\ref{table2} (IEEE-RTS-73 and IEEE-118), and we assess the computational complexity of each screening approach. We set the uncertainty parameter $\eta= 0.035$ for the two test systems. For the umbrella constraint screening algorithm with polyhedral uncertainty sets, all principal components are utilized to determine the corresponding uncertainty sets.

First, one of the main advantages of the proposed approach is that in terms of run time, it is significantly faster than Roald’s method (BA), as seen in Table~\ref{table3}. 
Note that we assess the computational time while we consider all $K = N^\star$ principal components, where $N^\star$ is the number of nodes with uncertain net demand such that $N^\star \leq N$. The conventional UCD algorithm has lower computational cost compared to the enhanced data-driven versions (D1-UCD and D2-UCD); however, with B-UCD a higher number of umbrella constraints are recorded. Table~\ref{table3} clearly shows that the solution time of the BA is directly proportional to the number of lines in the power system, while our approach is simply not. The UCD runtime is mainly affected by the final number of discovered umbrella network constraints and the system size.


Fig.~\ref{fig:04} illustrates the impact of the number of retained principal components on the screening time and the umbrella constraints. We see that when using DPUS $P_1$ while retaining only 10\% (for IEEE-RTS-73) and 40\% (for IEEE-118) of the most dominant principal components---Fig.~\ref{fig:04c} and Fig.~\ref{fig:04d}---identified the fewest umbrella constraints. On the other hand, $P_2$ needs respectively 60\% and 50\% of the principal components to match the umbrella constraint set counts of $P_1$.

\begin{table}
\captionsetup{font=footnotesize}
\caption{\sc{Number of Umbrella Constraints Identified by UCD Iteration for Medium-Size Networks}} 
\centering 
\footnotesize{
\scalebox{0.95}{
\begin{tabular}{c c c c c c c c c c c c c} 
\hline
& \multicolumn{10}{c}{Number of constraints found} & Total\\
\hline 
Iteration number & 1 & 2 & 3 & 4 & 5 & 6 & 7 & 8 & 9 & 10 \\ 
\hline 
& \multicolumn{10}{c}{B-UCD} & \\ 
\hline
RTS-73                  & 10    & 5    & 2     & 2   & 2  & 1    & 1    & 1     & --   & -- &  24\\[1.0ex]
IEEE-118                & 22   & 11   & 9      & 7   & 3  & 2    & 2    & 2     & 2    & -- &  60 \\[1.0ex] 
\hline
& \multicolumn{10}{c}{D1-UCD} & \\
\hline
RTS-73     & 8    & 5    & 2     & 2   & 1  & 1    & --    & --     & --  & -- &  19\\[1.0ex]
IEEE-118   & 22   & 9    & 7     & 4   & 3  & 3    & 2    & 2     & 1  & 1     &  53\\[1.0ex]
\hline 
\end{tabular}}}
\label{table4} 
\end{table}

\begin{table}
\captionsetup{font=footnotesize}
\caption{\sc{Redundant Constraint Screening  for Medium-Size Networks}} 
\centering 
\scalebox{0.62}{
\begin{tabular}{ccccc} 
\hline
  & D1-UCD  & BA (\% change) & D1-UCD & BA (\% change) \\ [1.0ex]
\hline 
& \multicolumn{2}{c}{IEEE RTS-73} & \multicolumn{2}{c}{IEEE-118}  \\ [1.0ex]
\hline  
Retained constraints (\%)                   & 7.9   &  10.0 (+26.6\%)    & 14.2      & 16.1 (+11.8\%)    \\[1ex]
UC compute time reduction (\%)  & 11.5  & 14.0 (+18.0\%)  & 16.9     & 19.4 (+12.9\%)  \\[1ex]
\hline 
\end{tabular}}
\label{table5} 
\end{table}

The number of umbrella line flow constraints identified in each iteration of the UCD algorithm is provided in Table~\ref{table4}. The B-UCD flavor of the algorithm identifies the umbrella line flow constraints of the networks IEEE-RTS-73 and IEEE-118 in eight, and nine iterations, respectively. On the other hand, D2-UCD converges in six iterations for IEEE-RTS-73, and in 10 iterations for the IEEE-118 test system. The results further indicate that the maximum number of line flow limits which may become active simultaneously under D2-UCD is at most eight and 22 for IEEE-RTS-73 and IEEE-118, respectively. This is the case because in the first iteration UCD found that the maximum number of intersecting constraints for each problem. We note also that D2-UCD retains fewer constraints in comparison to B-UCD. This is because the DPUS used here are tighter than the box sets defined by \eqref{eq:P2l}. \par 

The first important observation to be made is that the UCD approaches yield for the 1440 time periods in the test set the same optimal UC solutions and costs as those obtained with BA constraint screening---since all necessary constraints are retained by UCD. Consequently, the performance of these methods, which is summarized in Table~\ref{table5}, is assessed and compared in terms of the percentage of umbrella constraints and the computational burden of the reduced UC running time relative to the computational time required to solve the full UC problem. 
Our proposed approach outperforms Roald’s method in terms of network constraints removal for both IEEE-RTS-73 and IEEE-118. Therefore, the proposed approach provides UC solution time computational savings 18.0\% and 12.9\% lower than Roald's method for IEEE-RTS-73 and IEEE-118, respectively. Note that UC performance results using D1-UCD are identical to those for D2-UCD; this explains why the D2-UCD results have been omitted from Table~\ref{table3}.

\subsubsection{Large-Size Network}
To show the efficiency of the proposed approach for more realistic cases, this section compares the simulation results of the different methods for the \text{CASE500\_pserc} test system which has 500 buses and 733 lines. All system data are available from the IEEE PES PGLib-OPF v17.08 benchmark library \cite{2019PGLib}. To keep computational times within reasonable limits, the optimality gap is set to 1\% when solving UC problems. Line-based decomposition is used to partition all versions of UCD problems (B-UCD, D1-UCD, D2-UCD) into smaller sub-problems as explained in Appendix~A. These sub-problems are independent from each other and have been considered in a sequential manner to allow for a fair comparison with Roald's benchmark approach. To this end, we have partitioned arbitrarily the full line constraint set which consists of 1466 line constraints into 15 blocks; each subgroup $\mathcal{L_{\kappa}}$ contains 150 line flow constraints to be considered by each sub-problem \eqref{eq:P2adecomp}--\eqref{eq:P2fdecomp}
except for the last block which only has 66 line constraints to be examined.\footnote{This choice is purely arbitrary here. Other partitioning approaches are possible.} The training set contains 7200 time periods instances while the test set was reduced to 480 time periods. Here, we consider keeping $K=50$ principal components to formulate both $P_1$ and $P_2$. Results in Table~\ref{table6} illustrate that B-UCD's screening time is lower than $P_1$'s and $P_2$'s by an average of 58\%. This is a direct result of the fact that the number of constraints induced by DPUS is significantly larger than that of the box constraint set \eqref{eq:P2l}.

On the other hand, results found in Table~\ref{table7} show that the proposed method involves reductions in  the number of umbrella constraints and screening time compared to the benchmark, consequently, the reduced UC problem and computational time decrease. The screening time of our proposed method is 67.9\% faster than the benchmark approach. Also, imposing a polyhedral uncertainty set reduces the number of retained UC umbrella constraints and computational time for solving the UC by 25\% and 28\%, respectively.

\begin{table}
\captionsetup{font=footnotesize}
\caption{\sc{Redundant Constraint Screening Time for the Large-Size Network}} 
\centering 
\footnotesize{
\scalebox{0.8}{
\begin{tabular}{c c c c c} 
\hline
Method & B-UCD & D1-UCD & D2-UCD   \\ [0.5ex] 
\hline 
Total screening time (s)                   & 219.5  & 535.3  & 548.0          \\[1ex]
Average screening time per block (s)    & 14.6   & 35.9  & 36.5        \\[1ex]
\hline 
\end{tabular}}}
\label{table6} 
\end{table}

\begin{table}
\captionsetup{font=footnotesize}
\caption{\sc{Large-Size Network Computational Results}} 
\centering 
\footnotesize{
\scalebox{0.9}{
\begin{tabular}{c c c} 
\hline
Method & D1-UCD & BA (\% change)   \\ [0.5ex] 
\hline 
Retained constraints (\%)              & 6.34   & 8.45 (+25.0\%)      \\[1ex]
Screening time (s)                     & 535.3  & 1668.1 (+67.9\%)    \\[1ex]
UC compute time reduction (\%)              & 11.32  & 15.81  (+28.0\%)   \\[1ex]
\hline 
\end{tabular}}}
\label{table7} 
\end{table}

\subsection{Inactive Constraints Screening}
Finally, we experiment on evaluating the benefits of adding the production cost upper bound \eqref{eq:ED1} (along with \eqref{eq:ED2} and \eqref{eq:ED3}) to the UC constraints prior to solving its corresponding UCD. Here, we test for the addition of the cost upper bound with (i) box uncertainty sets (ED-UCD), (ii) data-driven polyhedral uncertainty sets, namely $P_1$ (ED+D1-UCD). In this study, the same data sets used in prior sections are deployed for DPUS computation and cost upper bound characterization.

For the IEEE-RTS-73 system, there are two constraints set up to provide a piecewise cost upper bound, while we use a single linear constraint for the IEEE-118 test system---see Table~\ref{tab:fit} for the production cost upper bound parameters used. Here these choices were made upon visual inspection of the historical production cost data plotted against observed system-wide net load for each system. As a first observation, the two economic-driven approaches provide for the test set the same optimal costs as obtained by the benchmark approach except for two instances (out of 1440) of the IEEE-RTS-73 which provides infeasible UC results. This happened as a result of two outliers from the historical production cost data points in the cost-net load regressions which were not covered by the upper bound constraints. The results in Table~\ref{table8} show that ED+D1-UCD leads to considerable reductions in both the number of retained constraints and computational burden in comparison to ED-UCD. For the IEEE-RTS-73 test system, imposing a data-driven uncertainty set  combined with the economic-driven constraints reduces retained constraints and computational time by 50\% and 32\%, respectively. For the IEEE-118 test system, ED+D1-UCD obtains reduced UC problems with 55\% fewer constraints than ED-UCD with a computational time lower by 43\%. In that case none of the 1440 cases tested led to infeasible UC results.

\begin{table}
    \centering
    \captionsetup{font=footnotesize}
    \caption{\sc{Parameters of Production Cost Upper Bounds}}
    \label{tab:fit}
    \footnotesize{
    \begin{tabular}{ccccc}
    \hline
    Test system & $a_0$ & $b_0$ & $\Delta$ & $\Gamma$ \\
    \hline
    IEEE-RTS-73 & $-6.03 \times 10^4$ & 24.64 & 5.0 & 0 \\
     & $-8.80 \times 10^4$ & 27.79 & & \\[0.5ex]
    IEEE-118 & $-1.35 \times 10^4$ & 18.9 & 3.7 & 0 \\[0.5ex]
    CASE500\_pserc & $-8.84 \times 10^4$ & 20.52 & 3.6 & 0\\
    & $-1.03 \times 10^5$ & 21.44 & & \\
    & $-2.63 \times 10^5$ & 30.00 & & \\
    \hline
    \end{tabular}}
\end{table}


For the \text{CASE500\_pserc} case, the cost-driven approach has been applied by fitting three piecewise linear segments; see Table~\ref{tab:fit}. The results reveal the lowest percentage of umbrella constraint with ED+D1-UCD slightly below 1\%, indicating that this network is not heavily congested for the patterns of net load it has to handle. This results in the modest UC computation time from ED+D1-UCD which is yet 50\% faster than that obtained after running ED-UCD  albeit for a slightly higher screening time. However, unlike in the redundant screening case, both approaches lead to few infeasible UC solutions (out of 500 instances) since some original umbrella UC constraints have been misclassified by overly tight production cost upper bounds. These can be addressed by further tuning the parameters $\Delta$ and $\Gamma$ to relax the cost upper bounds.

\begin{figure}
\begin{center}
\includegraphics[width=6.0cm]{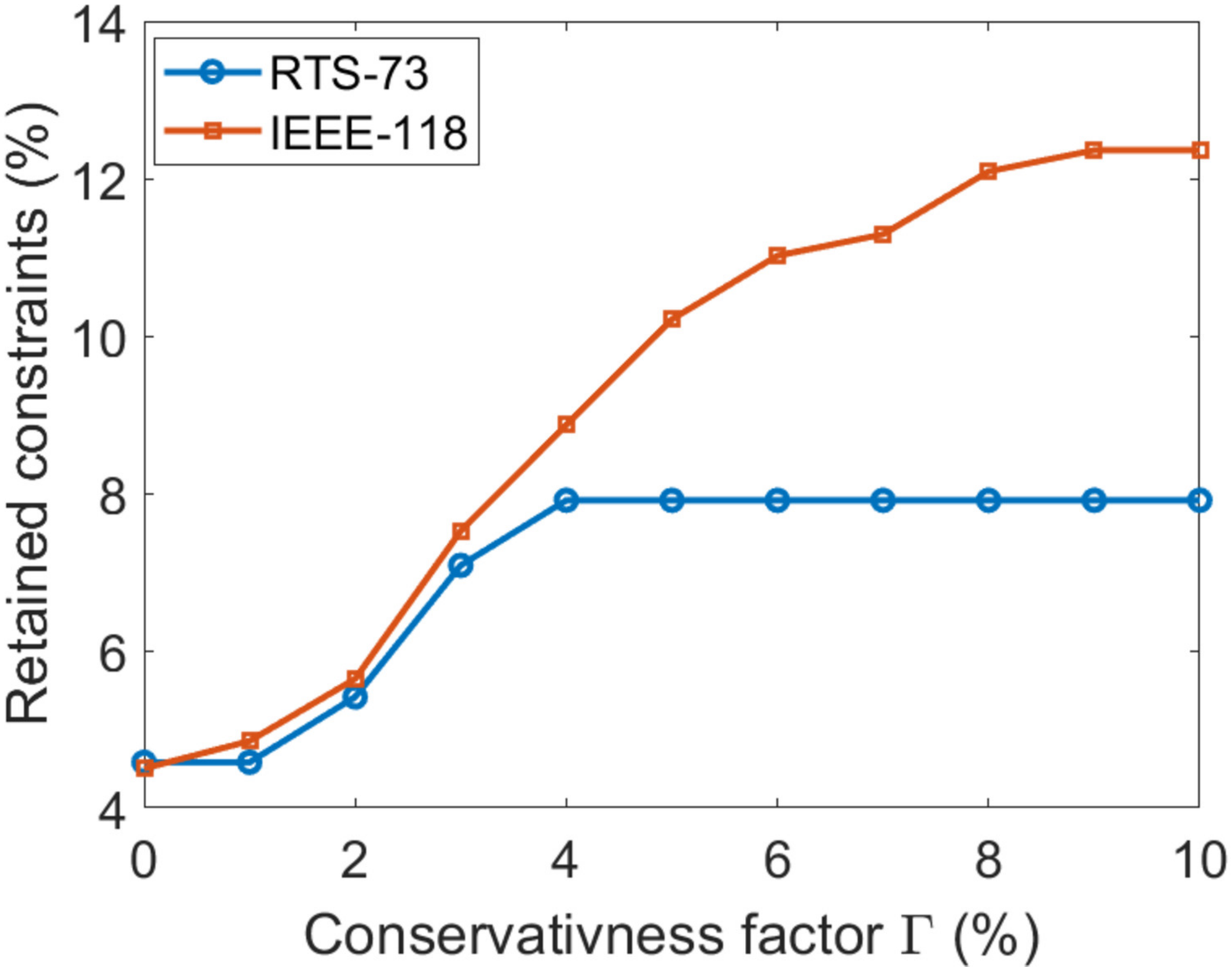}
\caption{Retained umbrella constraints in ED+D1-UCD as a function of $\Gamma$.}
\label{Fig4}
\end{center}
\end{figure}

To that end, Fig.~\ref{Fig4} shows how the conservativeness factor $\Gamma$ influences the constraint screening performance in ED+D1-UCD for the IEEE-RTS-73 and IEEE-118 cases. The conservativness factor increases the fitted slope values and yields higher numbers of retained umbrella constraints. For the IEEE-RTS-73 system, with conservativness factor values above 4\%, constraint screening results coincide with the number constraints found with D1-UCD as seen in Table~\ref{table4}. While for the IEEE-118 case, similar results are obtained when $\Gamma \geq 12\%$.

\begin{table}
\captionsetup{font=footnotesize}
\caption{\sc{Inactive Constraints Screening}} 
\centering 
\footnotesize{
\scalebox{0.8}{
\begin{tabular}{c c c} 
\hline
Method            & ED+D1 & ED (\% change)  \\ [0.5ex]
\hline 
 & \multicolumn{2}{c}{IEEE RTS-73}  \\ [0.5ex]
\hline  
Retained constraints (\%)              & 4.6  & 9.2 (+50.0\%)    \\[1ex]
Screening time (s)                         & 57.0   & 40.8 ($-$39.8\%) \\[1ex]
Number of infeasibilities                      & 2                  & 2      \\[1ex] 
UC compute time reduction (\%)             & 8.0 & 11.9 (+32.8\%)    \\[1ex]
\hline 
& \multicolumn{2}{c}{IEEE-118}  \\ [0.5ex]
\hline 
Retained constraints (\%)              & 4.5   & 10.2 (+55.2\%)          \\[1ex]
Screening time (s)                     & 19.7 & 40.3 (+51.1\%)       \\[1ex]
Number of infeasibilities                      & 0              & 0        \\[1ex] 
UC compute time reduction (\%)         & 6.8  & 12.0 (+43.7\%)    \\[1ex]
\hline 
& \multicolumn{2}{c}{\text{CASE500\_pserc}}  \\ [0.5ex]
\hline
Retained constraints (\%)              & 0.8    & 1.7 (+52.9\%)          \\[1ex]
Screening time (s)                     & 1517.4 & 1412.1 ($-$6.9\%)    \\[1ex]
Number of infeasibilities                      & 8      & 6                     \\[1ex]
UC compute time reduction (\%)         & 1.8  & 3.6 (+50.6\%)    \\[1ex]
\hline
\end{tabular}}}
\label{table8} 
\end{table}

\section{Conclusion}
A widespread observation in power system operation and planning optimization is that only a very small proportion of transmission line limitations are ever binding. In this paper, we proposed a data-driven umbrella constraint discovery approach that takes advantage of historical information to disregard redundant constraints in unit commitment problem formulations robust to uncertain net demand (demand less non-dispatchable generation). 

This approach requires first a characterization of polyhedral uncertainty sets using historical net demand observations and their past forecasts along with principal component analysis. We demonstrated how data-driven polyhedral uncertainty sets offer a good robustness compromise in comparison to cruder box-like uncertainty sets since they acknowledge cross-node correlations between demand and non-dispatchable generation.


Furthermore, we presented a proposal to refine the process of umbrella constraint discovery with the objective of predicting which of the umbrella constraints will be active in a unit commitment instance. With that prediction at hand, we demonstrate significant unit commitment solution speed-ups. This refinement is rendered possible by recognizing that the original set of umbrella constraints of a problem will intersect with a set of production cost upper bounds, whose computation is based on historical unit commitment solutions.


Future work is needed to adapt the proposed approach to the multi-period unit commitment problem, which involves taking into account inter-temporal constraints of generating units as well as with  security-constrained unit commitment problems. 

\appendices
\section{UCD Decomposition}
The computational burden of MIP in the UCD algorithm can be further reduced through the decomposition technique \cite{Amir2}, also called \emph{enhanced UCD (EUCD)}, as one can search for umbrella constraints in computationally-manageable
constraint blocks $\mathcal{L_{\kappa}}$ of the initial optimization $\cup_\kappa \mathcal{L}_\kappa=\{1, \ldots, L\}$ and $\mathcal{L}_\kappa \cap \mathcal{L}_{\kappa^{\prime}}=\emptyset$ for all $\kappa \neq \kappa^{\prime}$. To
perform this decomposition on UCD, one needs to consider the entire blocks of constraints \eqref{eq:P2b}--\eqref{eq:P2e} and \eqref{eq:P2k}--\eqref{eq:P2l} in addition to one of the blocks $\mathcal{L_{\kappa}}$ of the block of constraints in \eqref{eq:P2f}--\eqref{eq:P2i}. Furthermore, this step can be executed by parallel processing to expedite the solution of UCD problem. Hence, we solve \eqref{eq:P2adecomp}--\eqref{eq:P2fdecomp} for each $\mathcal{L_{\kappa}}$

\begin{equation}
\min \sum_{l^{\prime} \in \mathcal{L_{\kappa}}} (v^{+}_{l^{\prime}} + v^{-}_{l^{\prime}})\label{eq:P2adecomp}
\end{equation}

Subject to:
\begin{align}
\eqref{eq:P2b}-\eqref{eq:P2e} ,& \eqref{eq:P2k}-\eqref{eq:P2l} 
\label{eq:P2bdecomp}\\
\sum_{n=1}^N h_{l^{\prime} n} q_{n} + z^{+}_{l^{\prime}} & \geq f_{l^{\prime}}^{\max}, & \forall l^{\prime} \in \mathcal{L_{\kappa}}
\label{eq:P2cdecomp}\\
- \sum_{n=1}^N h_{l^{\prime} n} q_{n} + z^{-}_{l^{\prime}} & \geq f_{l^{\prime}}^{\max},& \forall l^{\prime} \in \mathcal{L_{\kappa}}
\label{eq:P2ccdecomp}\\
v_{l^{\prime}} - \frac{z^{+}_{l^{\prime}}}{\Omega} &\geq 0, & \forall l^{\prime} \in \mathcal{L_{\kappa}}
\label{eq:P2ddecomp}\\
v^{+}_{l^{\prime}} - \frac{z^{-}_{l^{\prime}}}{\Omega} &\geq 0, & \forall l^{\prime} \in \mathcal{L_{\kappa}}
\label{eq:P2ddecomp}\\
z^{+}_{l^{\prime}}, z^{-}_{l^{\prime}} &\geq 0,     & \forall l^{\prime} \in \mathcal{L_{\kappa}}
\label{eq:P2edecomp}\\
v^{+}_{l^{\prime}},v^{-}_{l^{\prime}} &\in\{0,1\}, & \forall l^{\prime} \in \mathcal{L_{\kappa}}
\label{eq:P2fdecomp}
\end{align}
In this case, the number of binary variables induced for examining line constraints per sub-problem is 2$\lvert\mathcal{L_{\kappa}}\lvert$, where $\lvert\mathcal{L_{\kappa}}\lvert$ is the cardinality of the subset $\mathcal{L_{\kappa}}$. By considering all the umbrella constraints identified in each block, we obtain all the umbrella constraints of the original UC problem.

\bibliography{References} 
\bibliographystyle{ieeetr}
\vspace{12pt}
\end{document}